\documentclass{amsart}

\usepackage{amssymb}
\usepackage{graphicx}
\usepackage{tikz}
\usepackage{tikz-cd}
\usepackage[latin1]{inputenc}

\tikzstyle{startstop}=[rectangle,rounded corners,draw=black,fill=green!20]
\tikzstyle{process}=[rectangle,draw=black,fill=green!20]
\tikzstyle{arrow}=[thick,->,>=latex]

\newtheorem{thm}{Theorem}[section]
\newtheorem{cor}[thm]{Corollary}
\newtheorem{lem}[thm]{Lemma}
\newtheorem{prop}[thm]{Proposition}

\theoremstyle{definition}
\newtheorem{dfn}[thm]{Definition}
\theoremstyle{remark}
\newtheorem{rem}[thm]{Remark}
\numberwithin{equation}{section}
\newcommand{\skcrlo}{\raisebox{-0.25\height}{\includegraphics[width=0.5cm]{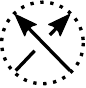}}}
\newcommand{\skcrro}{\raisebox{-0.25\height}{\includegraphics[width=0.5cm]{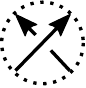}}}
\newcommand{\skcrvo}{\raisebox{-0.25\height}{\includegraphics[width=0.5cm]{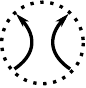}}}

\newcommand{\dL}{d\mspace{-2.5mu}l \mathbb{L}}
\begin{document}

\title[]{Braids for Knots in $S_{g} \times S^{1}$ and the affine Hecke algebra}%
\author{Seongjeong Kim}

\address{Seongjeong Kim, Jilin university}%
\email{kimseongjeong@jlu.edu.cn}%


\begin{abstract}
In \cite{Kim} it is shown that for an oriented surface $S_{g}$ of genus $g$ links in $S_{g} \times S^{1}$ can be presented by virtual diagrams with a decoration, called {\em double lines}. In this paper, first we define braids with double lines for links in $S_{g}\times S^{1}$. We denote the group of braids with double lines by $VB_{n}^{dl}$. The Alexander and Markov theorems for links in $S_{g}\times S^{1}$ can be proved analogously to the work in \cite{NegiPrabhakarKamada}. We show that if we restrict our interest to the group $B_{n}^{dl}$ generated by braids with double lines, but without virtual crossings, then the Hecke algebra of $B_{n}^{dl}$ is isomorphic to the affine Hecke algebra. Moreover, we define a Markov trace from the affine Hecke algebra to the Kauffman bracket skein module of $S^{2}\times S^{1}$.
\end{abstract}

\maketitle

\section{Introduction}

One of the generalizations of classical knot theory is {\em virtual knot theory}. It was discovered by Kauffman in 1996 \cite{Kauffman-virtual1,Kauffman-virtual2}. It was motivated by the study of knots in a thickened surface and a Gauss diagram. In particular, virtual knots can be presented by diagrams with two kinds of crossings, classical and virtual crossings, for example, see \cite{CarterKamadaSaito, KamadaKamada, Kauffman}. As a generalization of classical knot theory, various properties and applications of virtual knot theory have been found (cf. \cite{Manturov-Ilyutko_vir_book}).

On the other hand, braid theory connects knot theory with several research areas, for example, 3-manifolds theory and quantum field theory. And braid theory is connected with classical knot theory by Alexander theorem and Markov theorem. They state that any link can be presented by braids up to conjugation, stabilization and destabilization. Using this fact, the Jones polynomial was constructed via a Markov trace map on the Hecke algebra in 1987 \cite{Jones}. In \cite{GeckLambropoulou} by M. Geck and S. Lambropoulou Jones-type invariants for oriented links are constructed with trace maps on Iwahori-Hecke algebra of type B.

Analogously, virtual knots can be studied by using virtual braid theory. 
It is proved that any virtual link can be presented by braids up to conjugation, stabilization and right-left virtual exchange moves \cite{KauffmanLambropoulou-virbraid, Kamada-braid}. 

 We are interested in links in $S_{g}\times S^{1}$, where $S_{g}$ is an oriented surface of genus $g$. In \cite{Kim}, it is proved that links in $S_{g} \times S^{1}$ can be presented by virtual knot diagrams with decorations, called {\em double lines}. Moreover the isotopy of two knots in $S_{g} \times S^{1}$ can be presented by finitely many local moves. It is an extension of the work in \cite{DabkowskiMroczkowski}. That is, virtual knots are contained in the theory of knots in $S_{g}\times S^{1}$. In \cite{Kim-winding-parity}, the Kuperberg theorem for knots in $S_{g}\times S^{1}$ is proved. Moreover, a tool to study knots in $S_{g}\times S^{1}$, called {\em a winding parity}, is introduced and it is shown that there exist infinitely many knots in $S_{g}\times S^{1}$ of diagrams with only one or two crossings. In \cite{Kim-essential}, by using the winding parity it is proved that two virtual knot diagrams are equivalent modulo moves for knots in $S_{g}\times S^{1}$ if and only if they are equivalent as virtual knots, that is, modulo generalized Reidemeister move.

This paper contributes to braid theory for links in $S_{g}\times S^{1}$. In Section 2, we introduce the basic notions of links in $S_{g} \times S^{1}$. In Section 3, we define the group of virtual braids with double lines. We state the Alexander and Markov theorems for links in $S_{g}\times S^{1}$. In Section 4, we restrict our interest to the group of classical braids with double lines. We show that the Hecke algebra generated by the group of braids with double lines is isomorphic to the affine Iwahori-Hecke algebra. In Section 5, by the closure of classical braids with double lines, we construct a Markov trace into the Kauffman bracket skein module (algebra) of $S_{g}\times S^{1}$. In Appendix, we prove the Alexander and Markov theorem for links in $S_{g}\times S^{1}$.

\section{Link in $S_{g} \times S^{1}$ and its diagrams}

Let $S_{g}$ be an orientable surface of genus $g$. Let us define links in $S_{g} \times S^{1}$ analogously to virtual links by using underlying surfaces as follows:
\begin{dfn}
Let $S_{g}$ be an oriented surface of genus $g$. {\em A link $L$  in $S_{g} \times S^{1}$} is a pair of $S_{g} \times S^{1}$ and a smooth embedding of a disjoint union of $S^{1}$'s into $S_{g} \times S^{1}$. We denote it by $(L, S_{g} \times S^{1})$. Each image of $S^{1}$ in $S_{g} \times S^{1}$ is called {\em a component} of $L$. A link of one component is called {\em a knot in $S_{g} \times S^{1}$.}
\end{dfn}

\begin{dfn}
Let $(L, S_{g} \times S^{1})$ and $(L', S_{g'} \times S^{1})$ be two links in $S_{g} \times S^{1}$ and $S_{g'} \times S^{1}$. Two links $(L, S_{g} \times S^{1})$ and $(L', S_{g'} \times S^{1})$ are {\em equivalent} if $(L', S_{g'} \times S^{1})$ can be obtained from $(L, S_{g} \times S^{1})$ by isotopy and stabilization/destabilization of $S_{g} \times S^{1}$.
\end{dfn}
By the {\em destabilization for $(L, S_{g} \times S^{1})$} we mean the following:
let $C$ be a non-contractible circle on the surface $S_{g}$ such that there exists a torus $T$ in $S_{g} \times S^{1}$ homotopic to the torus $C \times S^{1}$ and not intersecting the link. Then the destabilization consists of cutting the manifold $S_{g} \times S^{1}$ along the torus $C \times S^{1}$ and gluing two newborn components of boundary with $D^{2} \times S^{1}$.
The {\em stabilization for $S_{g} \times S^{1}$} is the converse operation to the destabilization.

First let us construct diagrams for $(L, S_{g}\times S^{1})$ on the surface $S_{g}$ as follows: let $L$ be an (oriented) link in $S_{g} \times S^{1}$. Assume that an orientation is given on $S^{1}$. Suppose that $x_{0} \in S^{1}$ is a point such that $S_{g} \times \{x_{0}\} \cap L$ is a finite set of points with no transversal intersections. 
\begin{figure}[h]
\begin{center}
 \includegraphics[width = 12cm]{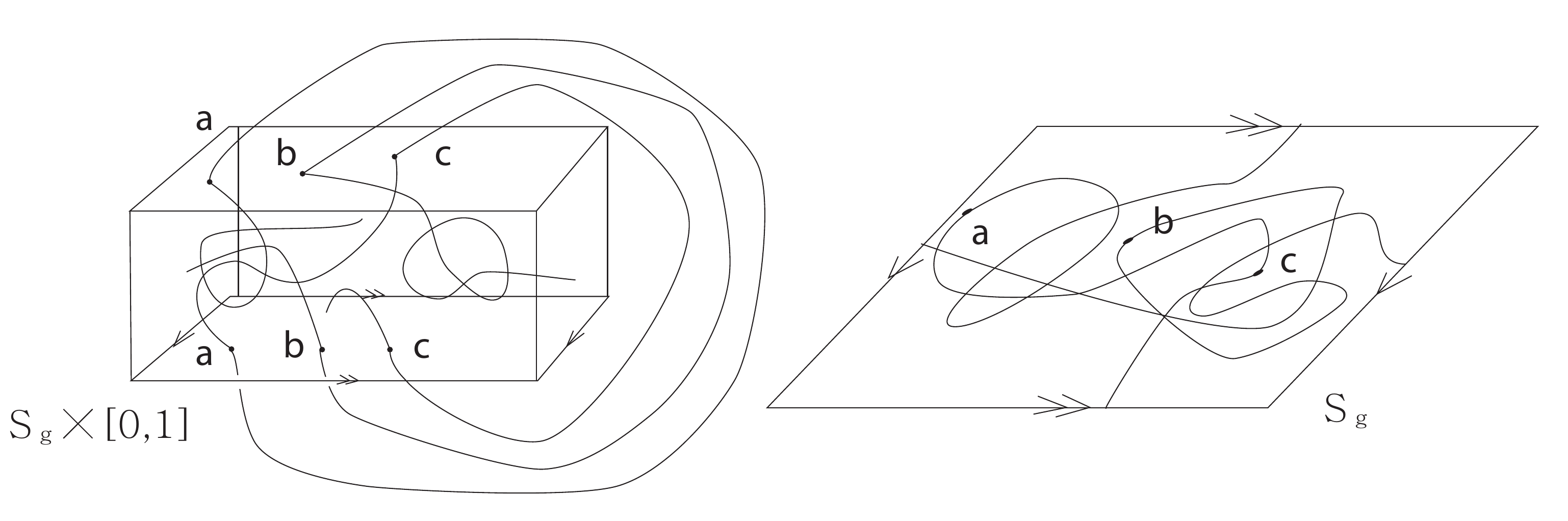}

\end{center}

\caption{Schematic figures of links in $S_{g} \times S^{1}$ and its projection on the surface $S_{g}\times \{0\}$}\label{fig:SS-diag-rel-1}
\end{figure}
Then there exists a natural diffeomorphism $f$\label{def:cutting-map} from $(S_{g} \times S^{1} )\backslash$  $(S_{g} \times \{x_{0}\})$ to $S_{g} \times (0,1) \subset S_{g} \times [0,1]$. Let $M_{L}$ $=$ $\overline{f((S_{g} \times S^{1}) - (S_{g} \times \{x_{0}\}))} \cong S_{g} \times [0,1]$. Then $\overline{f(L)}$ in $M_{L}$ consists of finitely many circles and arcs with exactly two boundaries on $S_{g} \times \{0\}$ and $S_{g} \times \{1\}$.

Let $D_{\overline{f(L)}}$ be the image of a projection of $\overline{f(L)}$ on the $S_{g} \times \{0\}$. Notice that boundary points of $\overline{f(L)}$ are paired so as to be projected to the same point on $S_{g}\times \{0\}$; we call such points {\em vertices}. It follows that the diagram $D_{\overline{f(L)}}$ of $L$ on $S_{g}$ has vertices corresponding to two boundary points on $S_{g} \times \{0\}$ and $S_{g} \times \{1\}$ of $\overline{f(L)}$ as described in the right of Fig.~\ref{fig:SS-diag-rel-1}.

\begin{figure}[h!]
\begin{center}
 \includegraphics[width = 8cm]{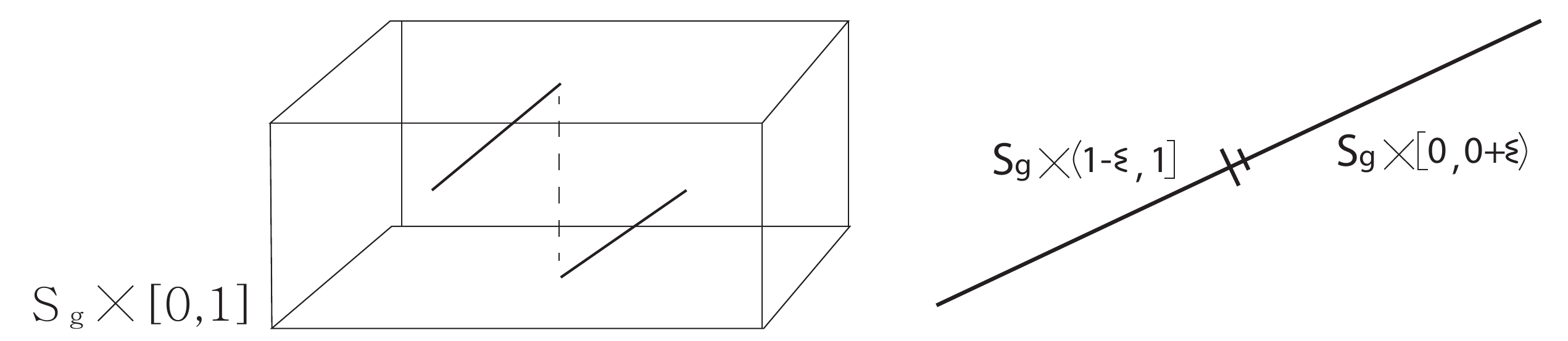}

\end{center}
\caption{Local image near to an intersection of a link with $S_{g} \times \{x_{0}\}$ and a corresponding double line}\label{fig:Vertex}
\end{figure}

Notice that two arcs near to a vertex are the images of arcs near to $S_{g} \times \{0\}$ and $S_{g} \times \{1\}$ in $M_{L} \cong S_{g} \times [0,1]$, respectively, as described in Fig.~\ref{fig:Vertex}. We change the point to two small lines so that if one of the lines is connected with an arc which is near to $S_{g} \times \{1\}$, then the line is longer than the another, as describe in Fig.~\ref{fig:Vertex}.

Since $D_{\overline{f(L)}}$ is a framed 4-valent graph with double lines on the surface $S_{g}$, which comes from $\overline{f(L)}$ in $S_{g} \times [0,1]$, one can give classical crossing information for each 4-valent vertex. That is, a link $L$ in $M_{L}$ has a knot diagram with double lines on $S_{g}$. We simply call it {\em a diagram on $S_{g}$ with double lines.} 

\begin{figure}[h]
\begin{center}
 \includegraphics[width = 12cm]{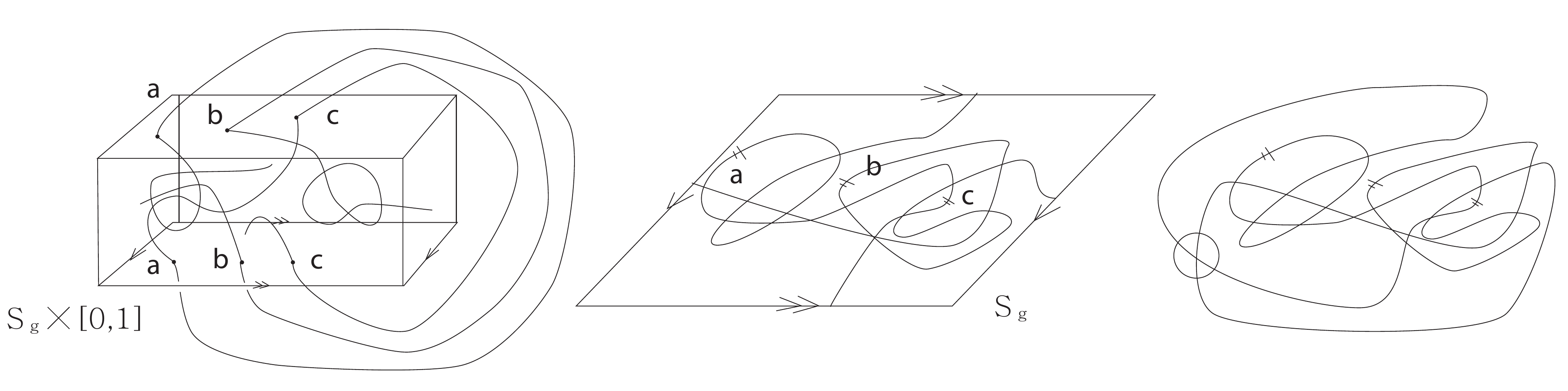}

\end{center}

\caption{Schematic figures of links in $S_{g} \times S^{1}$ and its projection on the plane}\label{fig:SS-diag-rel-2}
\end{figure}

Now, let us construct diagrams for links in $S_{g}\times S^{1}$ on the plane by using diagrams on $S_{g}$. For a link in $S_{g} \times S^{1}$ let a diagram $D$ on $S_{g}$ with double lines be given. We may assume that the diagram is drawn on $2g$-gon presentation of $S_{g}$ as in the middle of Fig.~\ref{fig:SS-diag-rel-2}. Connect points on boundaries of $2g$-gon with corresponding points by arcs outside $2g$-gon. By changing intersections between arcs outside $2g$-gon to virtual crossings, we obtain a diagram with double lines and virtual crossings, see the right in Fig.~\ref{fig:SS-diag-rel-2}. We call it {\it a diagram with double lines on the plane} for links in $S_{g}\times S^{1}$, or simply {\it a diagram for links in $S_{g}\times S^{1}$}.

The following theorem holds.
\begin{prop}[Kim (2018) \cite{Kim}\footnote{Diagrams on surface from $S_{g}\times S^{1}$ are studied in \cite{DabkowskiMroczkowski} to calculate Kauffman bracket skein module. They use the notion of ``arrows'' corresponding to double lines in this paper.}]\label{thm:diag_on_plane}
   Let $(L, S_{g} \times S^{1})$ and $(L', S_{g'} \times S^{1})$ be two links. Let $D_{L}$ and $D_{L'}$ be diagrams of $L$ and $L'$ on the plane, respectively. Then $L$ and $L'$ are equivalent if and only if $D_{L'}$ can be obtained from $D_{L}$ by applying finitely many moves in Fig.~\ref{fig:moves2}.
  \begin{figure}[h!]
\begin{center}
 \includegraphics[width = 12cm]{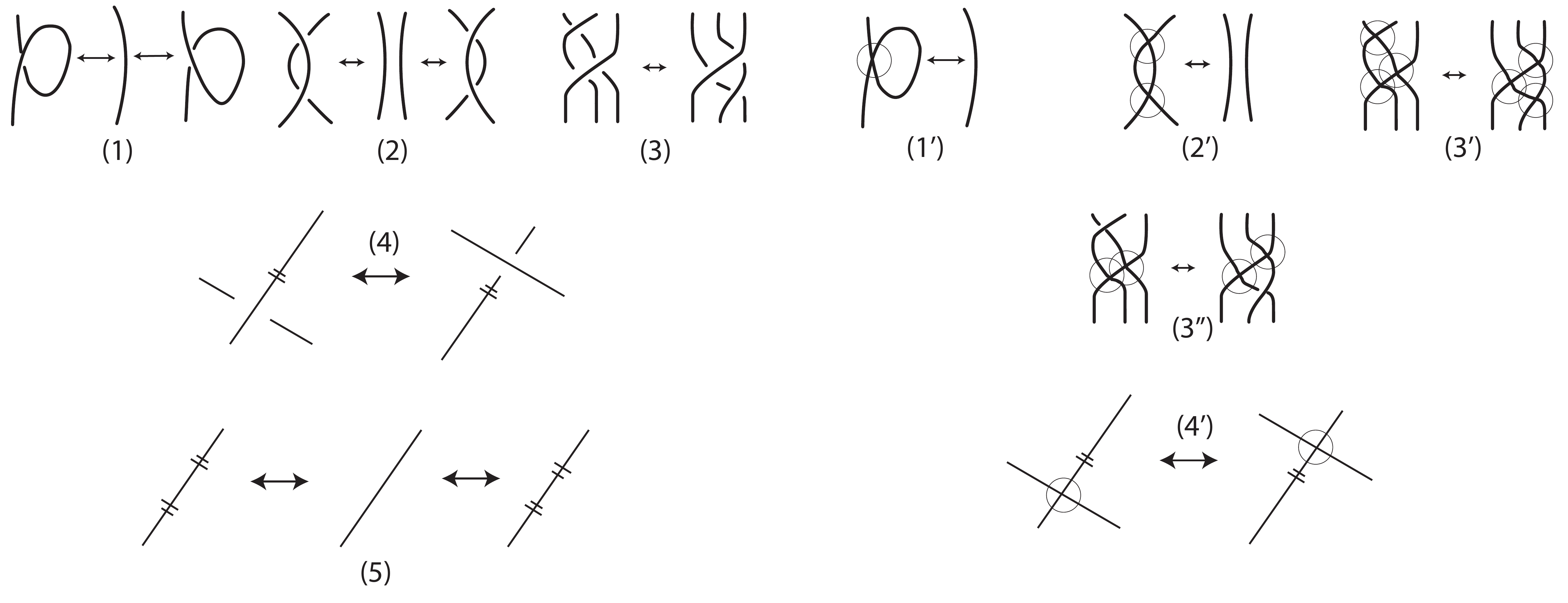}
\end{center}
 \caption{Moves for links in $S_{g}\times S^{1}$}\label{fig:moves2}
\end{figure}
\end{prop}

\begin{figure}[h]
\begin{center}
 \includegraphics[width = 9cm]{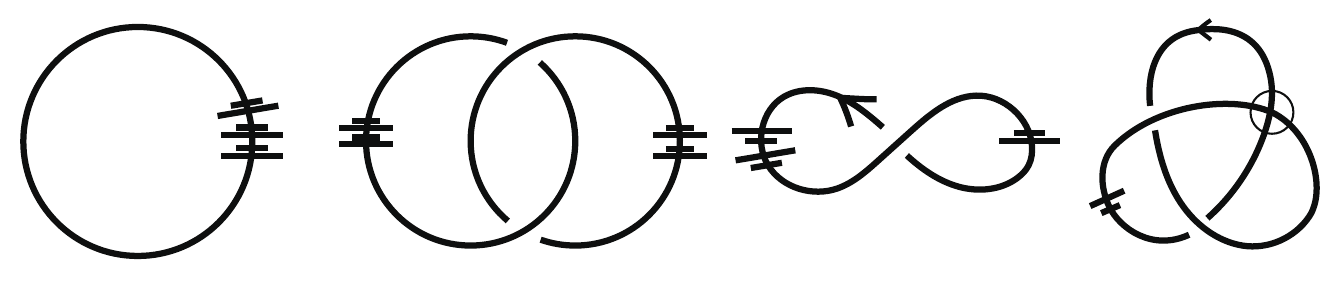}

\end{center}

\caption{Diagrams with double lines}\label{fig:exa-diag-dl}
\end{figure}

Simply speaking, knots and links in $S_{g}\times S^{1}$ can be presented by virtual knot diagrams with a double lines up to moves in Fig.~\ref{fig:moves2}, for example, see Fig.~\ref{fig:exa-diag-dl}.

\section{Virtual dl-braid group and its representation}

\subsection{Virtual braids with double lines}
\begin{dfn}
    {\em A virtual braid diagram $\beta$ with double lines on $n$ strands} is a union of smooth curves $b_{1}, \dots, b_{n}$ with double lines on $\mathbb{R} \times [0,1]$, called {\em strands} such that: 
    \begin{itemize}
        \item $\cup \partial(b_{i}) = \{x_{1}, \dots, x_{n}\} \times \{0,1\} \subset \mathbb{R} \times [0,1]$
        \item Each strand decreases monotonically from $\mathbb{R} \times \{1\}$ to $\mathbb{R} \times \{0\}$.
        \item Two strands intersect transversely. Each intersection is a classical crossing or a virtual crossing.
        \item There are only finitely many $t_{s} \in [0,1], s=1, \dots, m$ such that $\beta \cap \mathbb{R} \times \{t\}$ has a classical or virtual crossing or a double line. 
    \end{itemize}
For simplicity, we call such a braid {\em a virtual dl-braid diagram on $n$-strands}. We assume that each strand is oriented from $x_{s}\times \{1\}$ to $x_{t}\times \{0\}$, where $x_{s}\times \{1\}$ and $x_{t}\times \{0\}$ are the boundary of the strand.
\end{dfn}

\begin{dfn}
Two virtual dl-braid diagrams $\beta$ and $\beta'$ are equivalent if one can be obtained from another one by finite sequence of isotopies on $\mathbb{R} \times [0,1]$ preserving the conditions for dl-braids and the moves in Fig.~\ref{fig:dl-braid-moves}.
\begin{figure}[h!]
\begin{center}
 \includegraphics[width = 8cm]{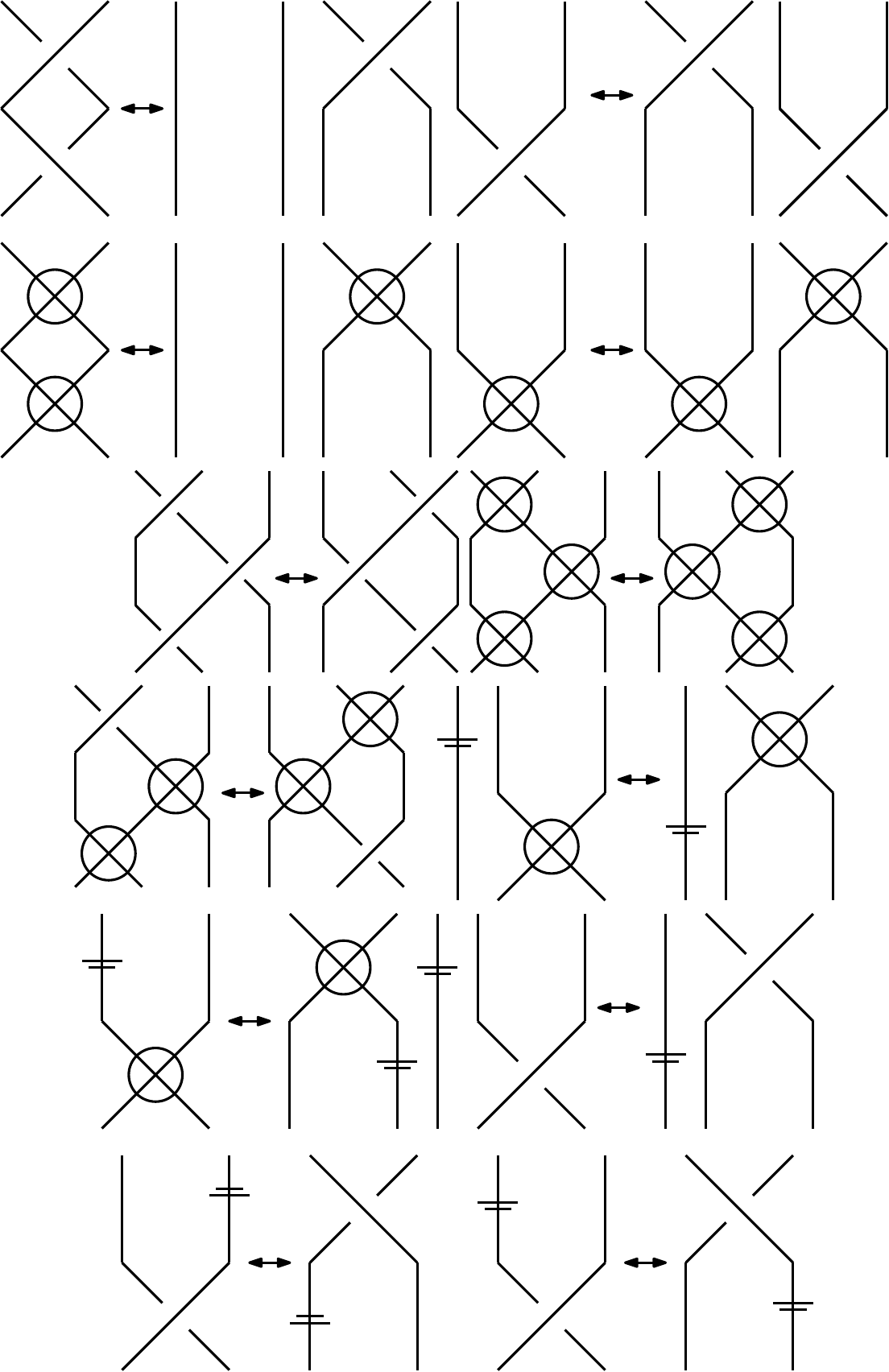}

\end{center}

 \caption{Moves for dl-braids}\label{fig:dl-braid-moves}
\end{figure}
\end{dfn}
Let $VB_{n}^{dl}$ be the set of equivalence classes of virtual dl-braids on $n$-strands. It has a natural binary operation $\beta_{1}\cdot \beta_{2}$ defined by placing $\beta_{1}$ on $\mathbb{R} \times [\frac{1}{2},1]$ and $\beta_{2}$ on $\mathbb{R} \times [0,\frac{1}{2}]$. One can show that $(VB_{n}^{dl},\cdot)$ is a group.

\begin{figure}[h!]
\begin{center}
 \includegraphics[width = 9cm]{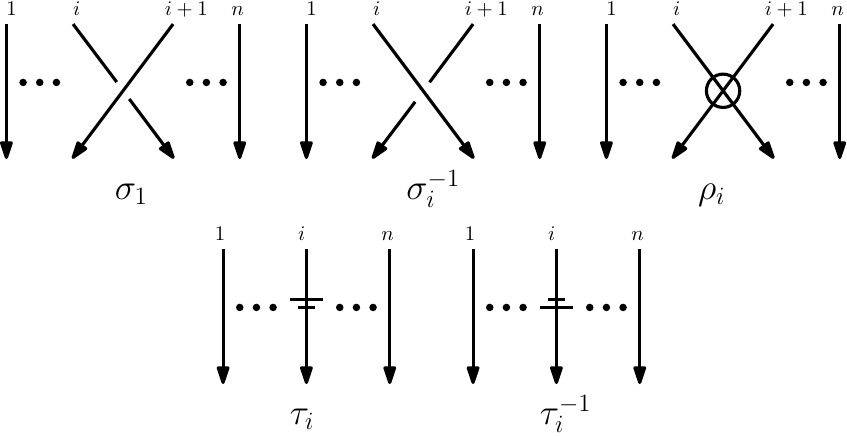}

\end{center}

 \caption{Generators for $VB_{n}^{dl}$}\label{fig:gen-dl-braid}
\end{figure}

Let us associate a positive (negative) crossing between $i$-th and $i+1$-th to $\sigma_{i}$ ($\sigma_{i}^{-1})$, a virtual crossing between $i$-th and $i+1$-th to $\rho_{i}$. And we associate a double line on $i$-th strand with positive (negative) sign to $\tau_{i}$ ($\tau_{i}^{-1}$), see Fig.\ref{fig:gen-dl-braid}.
Then one can easily see that $VB_{n}^{dl}$ can be defined by the following group presentation.

\begin{dfn}
    The $n$-strand dl-braid group $VB_{n}^{dl}$ is a group generated by generators $\{\sigma_{i}\}_{i=1}^{n-1}$, $\{\rho_{i}\}_{i=1}^{n-1}$ and $\{\tau_{i}\}_{i=1}^{n}$ and relations:
    \begin{enumerate}
        \item $\sigma_{i}\sigma_{j} = \sigma_{j}\sigma_{i}$, for $|i-j|\geq 2$.
        \item $\sigma_{i}\sigma_{i+1}\sigma_{i} = \sigma_{i+1}\sigma_{i}\sigma_{i+1},$ for $i=1,\dots, n-1$.
        \item $\rho_{i}^{2} =1$, for $i=1,\dots,n-1$.
        \item $\rho_{i}\rho_{j} = \rho_{j}\rho_{i}$, for $|i-j|\geq 2$.
        \item $\rho_{i}\rho_{i+1}\rho_{i} = \rho_{i+1}\rho_{i}\rho_{i+1},$ for $i=1,\dots, n-1$.
        \item $\sigma_{i}\rho_{i+1}\rho_{i} = \rho_{i+1}\rho_{i}\sigma_{i+1},$ for $i=1,\dots, n-1$.
        \item $\tau_{i}\tau_{j} = \tau_{j}\tau_{i}$ for $i,j=1,\dots, n$.
        \item $\sigma_{i}\tau_{j} = \tau_{j}\sigma_{i}$, if $j\neq i, i+1$.
        \item $\rho_{i}\tau_{j} = \tau_{j}\rho_{i}$ for $i=1,\dots n-1$, $j=1,\dots, n$.
        \item $\tau_{i}\sigma_{i} = \sigma_{i}^{-1}\tau_{i+1}$ for $i=1,\dots,n-1$.
    \end{enumerate}
\end{dfn}

The closure $\hat{\beta}$ of a virtual dl-braid diagram $\beta$ is defined in the standard way in knot theory. Then $\hat{\beta}$ is a diagram with double lines. It is clear that if $\beta$ and $\beta'$ are equivalent virtual dl-braids, then $\hat{\beta}$ and $\hat{\beta'}$ are equivalent diagrams with double lines. One important task with virtual dl-braid is to prove Alexander and Markov theorems for virtual dl-braids and diagrams with double lines. In \cite{Kamada-braid} Markov theorem for virtual braids is proved by S. Kamada. In \cite{NegiPrabhakarKamada} Alexander and Markov theorems for twisted virtual braids are proved by K. Negi, M. Prabhakar, S. Kamada. Since twisted virtual braids are virtual braids with a decoration, so called {\em bar}, it is similar to our situation. Since, roughly speaking, the detour move works for dl-braids, one can obtain Alexander and Markov theorems for dl-braids in similar way. We state them below. 

\begin{cor}
 Let $D$ be a diagram with double lines. There exists a virtual dl-braid $\beta$ such that $\hat{\beta}$ is equivalent to $D$.
\end{cor}

\begin{figure}[h!]
\begin{center}
 \includegraphics[width = 6cm]{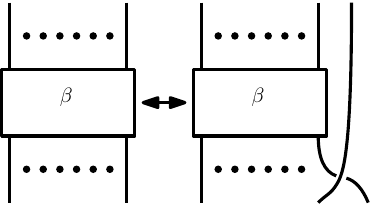}

\end{center}

 \caption{Stabilization and destabilization of a braid $\beta\sigma_{n} = \beta$ where $\beta \in B_{n}^{dl}$}\label{fig:stabilization}
\end{figure}

\begin{figure}[h!]
\begin{center}
 \includegraphics[width = 8cm]{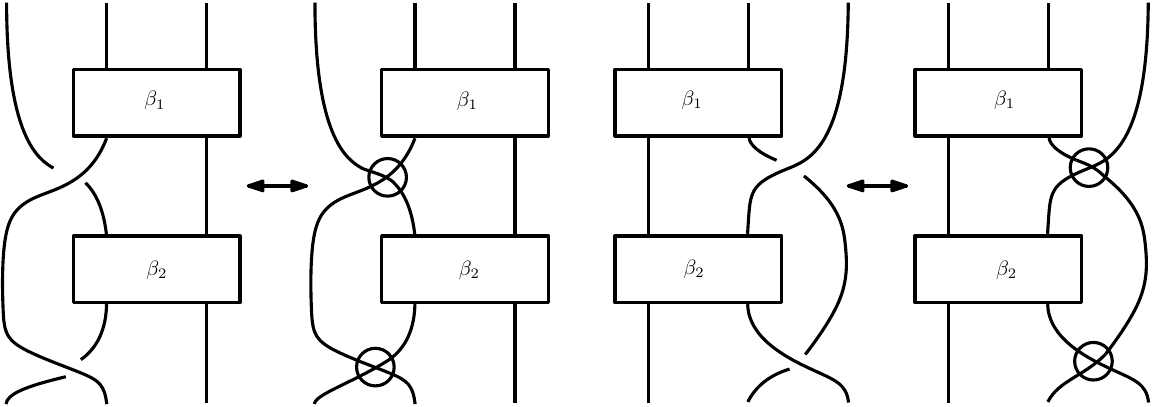}

\end{center}

 \caption{Left and right virtual exchange move: $\beta_{1}\sigma_{1}\beta_{2}\sigma_{1} = \beta_{1}\rho_{1}\beta_{2}\rho_{1}$ and $\beta_{1}\sigma_{n-1}\beta_{2}\sigma_{n-1} = \beta_{1}\rho_{n-1}\beta_{2}\rho_{n-1}$}\label{fig:vir-exchange}
\end{figure}

\begin{cor}
 Let $D, D'$ be diagrams with double lines. Let $\beta$ and $\beta'$ be two braids such that $D= \hat{\beta}$ and $D'=\hat{\beta'}$. $D$ and $D'$ are equivalent if and only if $\beta'$ can be obtained from $\beta$ by the following moves:
 \begin{itemize}
     \item Moves in Fig.~\ref{fig:dl-braid-moves}
     \item Conjugations
     \item Stabilization and destabilization depicted in Fig.~\ref{fig:stabilization}.
     \item Left and right virtual exchange moves depicted in Fig.~\ref{fig:vir-exchange}.
 \end{itemize}
\end{cor}

For completeness of the paper, the above two statements will be proved in Appendix.

\section{Hecke algebra generated classical dl-braid group and affine Hecke algebra}

\begin{dfn}\label{dfn:classical-dlbraid-gp}
 The $n$-strand classical dl-braid group $B_{n}^{dl}$ is a group generated by generators $\{\sigma_{i}\}_{i=1}^{n-1}$ and $\{\tau_{i}\}_{i=1}^{n}$ and relations:
    \begin{enumerate}
        \item $\sigma_{i}\sigma_{j} = \sigma_{j}\sigma_{i}$, for $|i-j|\geq 2$.
        \item $\sigma_{i}\sigma_{i+1}\sigma_{i} = \sigma_{i+1}\sigma_{i}\sigma_{i+1},$ for $i=1,\dots, n-1$.
        \item $\tau_{i}\tau_{j} = \tau_{j}\tau_{i}$ for $i,j=1,\dots, n$.
        \item $\sigma_{i}\tau_{j} = \tau_{j}\sigma_{i}$ for $j \neq i,i+1$.
        \item $\tau_{i}\sigma_{i} = \sigma_{i}^{-1}\tau_{i+1}$ for $i=1,\dots,n-1$.
    \end{enumerate} 
\end{dfn}

\begin{rem}
    Note that in this paper, we do not insist that $B_{n}^{dl}$ is a subgroup of $VB_{n}^{dl}$. However, the answer is positive, since it is known that the classical braid group is a subgroup of the virtual braid group. It is proved by using parity and it is not affected by double lines and over/under information of crossings. We will study this in a separate paper.
\end{rem}

\begin{dfn}[\cite{Schiffmann}]\label{dfn:Affine-Hecke}
    The Iwahori-Hecke algebra $\hat{\mathcal{H}}_{n}$ associated to $\hat{\mathcal{G}}_{r} = \mathcal{G}_{r} \ltimes \mathbb{Z}^{r}$, where $\mathcal{G}_{r}$ is the symmetric group on $r$ elements, is the unital $\mathbb{C}[v^{\pm 1}]$-algebra generated by $T_{i}^{\pm 1}$, $X_{j}^{\pm1}$ where $i =1,\dots,n-1$ and $j=1,\dots,n$ with the following relations:
    \begin{enumerate}
        \item $T_{i}T_{i}^{-1} = 1 = T_{i}^{-1}T_{i}$.
        \item $T_{i}T_{j} = T_{j}T_{i}$ for $|i-j|\geq 2$,
        \item $T_{i}T_{i+1}T_{i} = T_{i+1}T_{i}T_{i+1}$
        \item $X_{i}X_{i}^{-1} =1 = X_{i}^{-1}X_{i}$,
        \item $X_{i}X_{j} = X_{j}X_{i}$ for $|i-j|\geq 2$,
        \item $T_{i}X_{j} = X_{j}T_{i}$ if $j\neq i, i+1$
        \item $T_{i}X_{i}T_{i} = v^{-2}X_{i+1}$,
        \item $(T_{i}+1)(T_{i}-v^{-2})=0.$
    \end{enumerate}
    It is also called {\em the affine Hecke algebra}
\end{dfn}

Let $\mathcal{H}[B_{n}^{dl}]$ be the Hecke algebra of $B_{n}^{dl}$, that is, it is the unital $\mathbb{C}[q^{\pm 
\frac{1}{2}}]$-algebra generated by $B_{n}^{dl}$ with the relation
\begin{equation}\label{eq:Hecke-rel}
    (\sigma_{i}-q^{\frac{1}{2}})(\sigma+q^{-\frac{1}{2}})=0.
\end{equation} 

\begin{thm}
    $\mathcal{H}[B_{n}^{dl}]$ is isomorphic to the affine Hecke algebra $\hat{\mathcal{H}}_{n}$.
\end{thm}

\begin{proof}
    Let us define a map $\phi$ from $\mathcal{H}[B_{n}^{dl}]$ to $\hat{\mathcal{H}}_{n}$ by $\phi(q^{\frac{1}{2}}\beta) = v^{-1}\phi(\beta)$, $\phi(\sigma_{i}) = vT_{i}$ and $\phi(\tau_{j}) = X_{j}$. Then it is clear that relation (1)-(4) in Definition~\ref{dfn:classical-dlbraid-gp} are preserved. For the relation (5) in Definition~\ref{dfn:classical-dlbraid-gp}, one can rewrite is as 
    $$\sigma_{i}\tau_{i}\sigma_{i}= \tau_{i+1}.$$
    Then 
    \begin{eqnarray*}
    \phi(\sigma_{i}\tau_{i}\sigma_{i}) &=& \phi(\sigma_{i})\phi(\tau_{i})\phi(\sigma_{i})\\
        &=&vT_{i}X_{i}(vT_{i}) \\
        &=& vT_{i}X_{i}T_{i}\\
        &=& v^{2}(v^{-2}X_{i+1})\\
        &=& X_{i+1} = \phi(\tau_{i+1}).
    \end{eqnarray*}
    For the equation (\ref{eq:Hecke-rel}),
    
    \begin{eqnarray*}
    \phi((\sigma_{i}-q^{\frac{1}{2}})(\sigma+q^{-\frac{1}{2}})) &=& (vT_{i}-v^{-1})(vT_{i}+v)\\
        &=&v^{2}(T_{i}-v^{-2})(T_{i}+1) \\
        &=& 0.
    \end{eqnarray*}
    Therefore $\phi$ is well-defined. Conversely, if we define a map by associating $v \mapsto q^{-\frac{1}{2}}$, $T_{i} \mapsto q^{\frac{1}{2}}\sigma_{i}$, $X_{j} \mapsto \tau_{j}$, then it is an inverse map of $\phi$ and hence it completes the proof.
\end{proof}

Therefore, we have a geometric description for affine Hecke algebra $\hat{\mathcal{H}}_{n}$ by using dl-braids. It means that there exists one-to-one correspondence between representations of $\hat{\mathcal{H}}_{n}$ and $\mathcal{H}[B_{n}^{dl}]$. In particular, one can find a trace map for $\hat{\mathcal{H}}_{n}$ by using the closure of braids.

\section{Skein module and a trace map on $\mathcal{H}[B_{n}^{dl}]$}

\begin{dfn}\label{kbsmdef}

Let $M$ be an oriented $3$-manifold, $R$ a commutative ring with unity, and $q^{\frac{1}{4}} \in R$ a fixed invertible element. Consider  the set of ambient isotopy classes of unoriented framed links (including the empty link $\varnothing$) in $M$, which we denote by $\mathcal{L}^{\mathit{fr}}$,  and the free $R$-module with basis $\mathcal{L}^{\mathit{fr}}$, denoted by $R\mathcal{L}^{\mathit{fr}}$.
Let $S_{2, \infty}^{\mathit{sub}}$ be the submodule of $R\mathcal{L}^{\mathit{fr}}$ generated by the following expressions: 

\begin{enumerate}

    \item the Kauffman bracket skein expression: $L_+ - q^{\frac{1}{4}}L_0 - q^{-\frac{1}{4}}L_{\infty}$,
    
    \item the trivial component expression: $L \sqcup {\pmb \bigcirc}  + (q^{\frac{1}{2}} + q^{-\frac{1}{2}})L$,
    
\end{enumerate}
\noindent
where $\pmb{\bigcirc}$ denotes the trivial framed knot in $M$ and the skein triple $(L_+$, $L_0$, $L_{\infty})$ denotes three framed links in $M$, which are identical except in a small $3$-ball in $M$ where they differ as illustrated in Figure \ref{skeintriple}.

\begin{figure}[ht]
    \centering
\[  \begin{minipage}{1 in} \includegraphics[width=\textwidth]{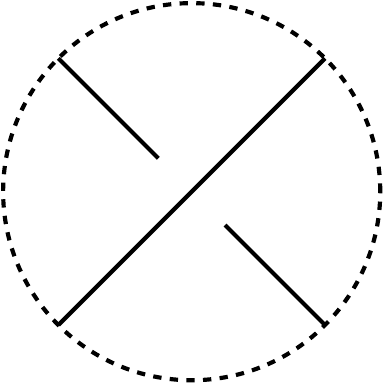} \vspace{-15pt} \[L_+\] \end{minipage} 
               \qquad
        \begin{minipage}{1 in}\includegraphics[width=\textwidth]{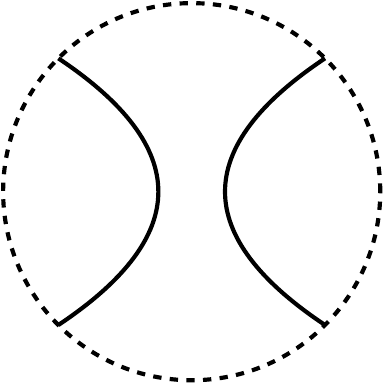} \vspace{-15pt} \[L_0\] \end{minipage}
         \qquad
        \begin{minipage}{1 in}\includegraphics[width=\textwidth]{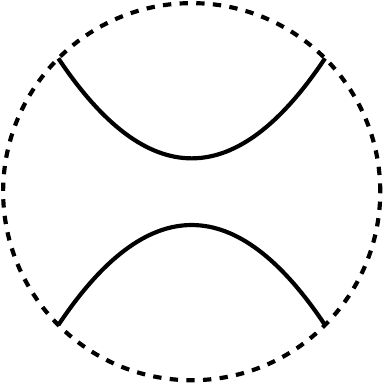} \vspace{-15pt} \[ L_\infty\]\end{minipage} 
        \]
\caption{Skein triple for the Kauffman bracket skein module.}
          \label{skeintriple}
        \end{figure}

The {\bf Kauffman bracket skein module} of $M$ is defined as the quotient: $$\mathcal{S}_{2,\infty}(M;R,A) = \frac{R\mathcal{L}^{\mathit{fr}}}{S_{2, \infty}^{sub}}.$$

\end{dfn}

In particular, when $M = S^{2}\times S^{1}$, we know that links in $S^{2} \times S^{1}$ can be presented by diagrams on $S^{2}$ with double lines up to local moves. Therefore, for $M = S^{2}\times S^{1}$, when we define skein module of $M$, `isotopy classes of links in $M$' can be replaced by `diagrams with double lines on $S^{2}$ modulo moves for them'. 

Now one can define a multiplication $L_{1}\cdot L_{2}$ in $S_{2,\infty}(S^{2}\times S^{1})$ by placing $L_{1}$ and $L_{2}$ disjointly on $S^{2}$. This multiplication makes $S_{2,\infty}(S^{2}\times S^{1})$ into a commutative algebra.

\begin{thm}
    Let us define $tr_{n} : \mathcal{H}[B_{n}^{dl}] \rightarrow S_{2,\infty}(S^{2}\times S^{1})$ by $$tr_{n}(\beta) = (-q^{\frac{1}{4}})^{-wr(\beta)}[\hat{\beta}],$$ 
    where $wr(\beta)$ is the writhe of the braid diagram. Then it satisfies 
    \begin{itemize}
        \item $tr_{n}(xy)= tr_{n}(yx)$
        \item $tr_{n+1}(\beta\sigma_{n}^{\pm1}) = -q\cdot tr_{n}(\beta)$ for $\beta \in B_{n}^{dl} \subset B_{n+1}^{dl}$.
    \end{itemize}
\end{thm}

\begin{proof}
    To show that $tr_{n}$ is well defined for each $n$, it suffices to show that the equation (\ref{eq:Hecke-rel}) is preserved. From the equation \ref{eq:Hecke-rel} one obtain 
    \begin{equation*}
        \sigma_{i}^{2}+(-q^{\frac{1}{2}}+q^{-\frac{1}{2}})\sigma_{i}-1=0.
    \end{equation*}
    Then  
    \begin{eqnarray*}
        tr_{n}((LHS)) &=& (-q^{\frac{1}{4}})^{2}\hat{\sigma_{i}^{2}}-(-q^{\frac{1}{4}})^{0}1+(-q^{\frac{1}{4}})^{1}(-q^{\frac{1}{2}}+q^{-\frac{1}{2}})\hat{\sigma_{i}}\\
        &=&-(-q^{\frac{1}{4}})(q^{\frac{1}{4}}\hat{\sigma_{i}^{2}}-q^{-\frac{1}{4}}1+(-q^{\frac{1}{2}}+q^{-\frac{1}{2}})\hat{\sigma_{i}}).
    \end{eqnarray*}
    Note that from the Kauffman bracket skein expression it follows that
    $$q^{\frac{1}{4}}\skcrro-q^{-\frac{1}{4}}\skcrlo=(q^{\frac{1}{2}}-q^{-\frac{1}{2}})\skcrvo.$$ Therefore, 
    $$tr_{n}((LHS)) = 0 =tr_{n}((RHS)).$$
    Since we close the braid, it follows that $tr_{n}(xy)= tr_{n}(yx)$. For the second property
    \begin{eqnarray*}
        tr_{n+1}(\beta\sigma_{n}) &=& (-q^{\frac{1}{4}})^{wr(\beta)+1}\hat{\beta\sigma_{n}} \\
        &=& (-q^{\frac{1}{4}})^{wr(\beta)+1}(-q^{\frac{1}{4}})^{3}\hat{\beta}\\
        &=&(-q)(-q^{\frac{1}{4}})^{wr(\beta)}\hat{\beta}\\
        &=&(-q)tr_{n}(\beta).
    \end{eqnarray*}
    
    
\end{proof}
That is, there is a unnormalized Markov-trace for affine Hecke algebra valued in the Kauffman bracket skein module of $S^{2}\times S^{1}$.

About the Kauffman bracket skein module of $S^{2}\times S^{1}$, the following fact is well-known.
\begin{prop}[\cite{HostePrzytycki}]
    $$\mathcal S_{2,\infty}(S^{2}\times S^{1})=\mathbb Z[q^{\pm \frac{1}{4}}]\oplus \bigoplus_{i=1}^{\infty} \mathbb Z[q^{\pm \frac{1}{4}}]/(1-q^{\frac{2i+4}{4}}),$$
\end{prop}

\begin{figure}[h]
\begin{center}
 \includegraphics[width = 6cm]{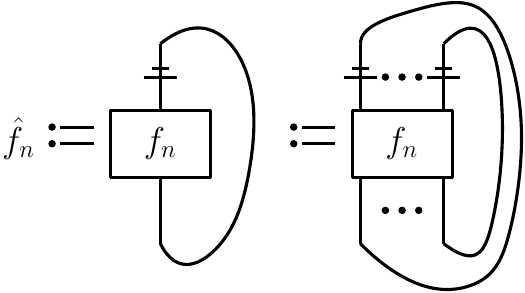}

\end{center}

\caption{$\hat{f}_{n}$ where $f_{n}$ is the Jones-Wenzl idempotent}\label{fig:Jones-Wenzl}
\end{figure}

In \cite{HostePrzytycki} generators for this decomposition are given by using a recurrence relation. In particular, when the generators are constructed, the Chebyshev polynomial with variable $x$, which represents $ \{*\}\times S^{1} $ for some point $* \in S^{2}$, plays an important role. In our expression, $x$ is $\hat{\tau}_{1}$, which is a trivial diagram with one double line. It follows that a value of $tr_{n}$ has a linear combination of diagrams $\hat{f_{n}}$ described in Fig.~\ref{fig:Jones-Wenzl}.

\begin{rem}\label{rem:trace-skein}
    In similar way, one can show that there exists a trace map from $tr_{n} : \mathcal{H}[B_{n}^{dl}] \rightarrow S_{Homflypt}(S^{2}\times S^{1}),$ where $S_{Homflypt}(S^{2}\times S^{1})$ is the Homflypt skein module with property  
    \begin{itemize}
        \item $tr_{n}(xy)= tr_{n}(yx)$
        \item $tr_{n+1}(\beta\sigma_{n}^{\pm1}) = tr_{n}(\beta)$ for $\beta \in B_{n}^{dl} \subset B_{n+1}^{dl}$.
    \end{itemize}
    Theoretically, we can define trace maps of affine Hecke algebra valued in any kinds of skein modules.
\end{rem}

\section{Further research}

The first question is to find and to classify trace maps. As mentioned in Remark~\ref{rem:trace-skein}, one can construct trace maps of affine Hecke algebra valued in skein modules. It helps us to understand the structure of affine Hecke algebra.
On the other hand, another natural question from the relationship between $\mathcal{H}[B_{n}^{dl}]$ and $\hat{\mathcal{H}}_{n}$ is to find a kind of affine Hecke algebra isomorphic to $\mathcal{H}[VB_{n}^{dl}]$. It must be related with the question how virtual braid group is related with Coxeter group. It is worth noting that the classical braid group is a subgroup of the virtual braid group. It means, one can expect that $B_{n}^{dl}$ is a subgroup of $VB_{n}^{dl}$ and hence $\hat{\mathcal{H}}_{n}$ is a subalgebra of $\mathcal{H}[VB_{n}^{dl}]$. We will discuss about this in a separate paper.

Related with the structure of $\mathcal{H}[VB_{n}^{dl}]$, one can expect that there exists a trace map from $\mathcal{H}[VB_{n}^{dl}]$ to KBSM of $S_{g}\times S^{1}$. The answer is negative, because stabilization and destabilization of underlying surface is required for $VB_{n}^{dl}$. It might be possible to consider the group of braids with double lines placed on a fixed oriented surface with a boundary. Then it might be related with affine Hecke algebra and gives rise to representations of affine Hecke algebra to KBSM of $S_{g}\times S^{1}$.

\section{Appendix : Alexander and Markov theorems for links in $S_{g}\times S^{1}$}
In this section, we prove the Alexander and Markov theorems for links in $S_g \times S^1$. 
The proof closely follows that of \cite{NegiPrabhakarKamada} for twisted braids and links, and we adopt the notation used there.

For a virtual diagram $D$ with double lines, 
\begin{itemize}
    \item Let $V(D)$ be the set of all classical crossings of $D$.
    \item Let $S(D)$ be the map from $V_{R}(D)$ to the set $\{\pm 1\}$ assigning the signs to classical crossings.
    \item Let $\dL(D)$ be the set of double lines on $D$.
    \item Let $N(v)$ be a regular neighborhood of $v$, where $v\in V(D) \cup \dL(D)$.
    \item For $c \in V(D)$ we denote by $c^{(1)},c^{(2)},c^{(3)}, c^{(4)}$ the four points of $\partial N(c)\cap D$ as described in Fig.~\ref{fig:bdy-of-crossing}.
    \begin{figure}[h]
\begin{center}
 \includegraphics[width = 6cm]{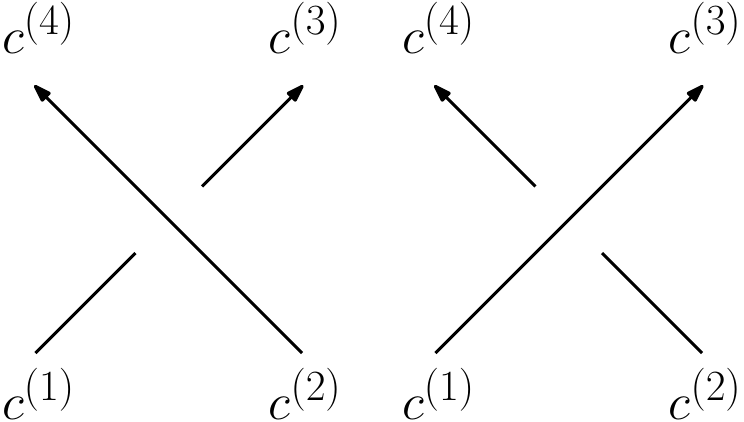}

\end{center}

\caption{Boudary points of a crossing}\label{fig:bdy-of-crossing}
\end{figure}
    \item For $d\in \dL(D)$ we denote by $d^{(1)}$ and $d^{(2)}$ the two points of $\partial N(d)\cap D$ as described in Fig.~\ref{fig:bdy-of-dl}.
    \begin{figure}[h]
\begin{center}
 \includegraphics[width = 3cm]{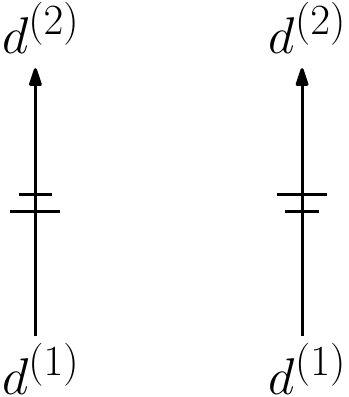}

\end{center}

\caption{Boudary points of a double line}\label{fig:bdy-of-dl}
\end{figure}
    \item Put $W(D) = cl(\mathbb{R}^{2}- \cup_{v \in V(D)\cup \dL(D)}N(v))$, where $cl$ is the closure.
    \item Let $V^{\partial}(D)= \{c^{(j)}~|~c\in V(D), 1\leq j\leq 4\}$, and $\dL^{\partial}(D) = \{d^{(j)}~|~d\in \dL(D), j= 1,2\}$.
    \item Let $D|_{W(D)}$ be the restriction of $D$ to $W(D)$, which is a union of oriented arcs and loops generically immersed in $W(D)$ such that the double points are virtual crossings of $D$, and the set of boundary points is $V^{\partial}(D) \cup \dL^{\partial}(D)$.
    \item Let $\mu(K)$ be the number of components of $K$.
    \item Define a subset $E(D)$ of $(V(D)^{\partial}\cup\dL^{\partial}(D))^{2}$ such that $(a,b) \in E(D)$ if and only if $K|_{W}$ has an oriented arc starting from $a$ and ends of $b$.

\end{itemize} 

{\em The Gauss data} of a diagram $D$ with double lines is the quintuple
$$(V(D),S(D),\dL(D),E(D),\mu(D)).$$

We say that two diagrams $D$ and $D'$ with double lines have {\em the same Gauss data} if $\mu(D)=\mu(D')$ and there exists a bijection $g : V(D)\cup \dL(D) \rightarrow V(D')\cup \dL(D')$ satisfying the following conditions:
\begin{itemize}
    \item $g(V(D))= V(D')$ and $g(\dL(D))= \dL(D')$.
    \item $g$ preserves the signs of classical crossings; $S(D)(c) = S(D')(g(c))$
    \item $(a,b) \in E(D)$ if and only if $(g^{\partial}(a),g^{\partial}(b)) \in E(D')$, where $g^{\partial} : V^{\partial}(D) \cup \dL^{\partial}(D) \rightarrow V^{\partial}(D') \cup \dL^{\partial}(D')$ is the bijection induced from $g$.
\end{itemize}
Suppose that $D'$ is a diagram with double lines with the same Gauss data with $D$. Then by an isotopy of $\mathbb{R}^{2}$ we can move $D'$ such that
\begin{itemize}
    \item $D$ and $D'$ are identical in $N(v)$ for every $v \in V(D)\cup \dL(D)$.
    \item $D'$ has no classical crossings and double lines in $W(D)$,
    \item there is a bijection between the arcs/loops of $D|_{W(D)}$ and those of $D'|_{W(D')}$ with respect to the endpoints of the arcs.
\end{itemize}
In this situation, we say that $D'$ is obtained from $D$ by replacing $D|_{W(D)}$.

\begin{lem}
    Let $D$ and $D'$ be diagrams with double lines. Let $W= W(D)= cl(\mathbb{R}^{2}- \cup_{v \in V(D)\cup \dL(D)}N(v))$.
    \begin{enumerate}
        \item If $D'$ is obtained from $D$ by replacing $D|_{W}$, then they are related by a finite sequence of isotopies of $\mathbb{R}^{2}$ with support $W$ and moves (1'),(2'),(3'), (3''), (4') in Fig. \ref{fig:moves2}.
        \item If two diagrams $D$ and $D'$ with double lines have the same Gauss data, then $D$ is equivalent to $D'$.
    \end{enumerate}
\end{lem}
\begin{proof}
    The proof is almost same to the proof of Lemma 3.3 in \cite{NegiPrabhakarKamada}. The idea is that a replacement of $D_{W}$ can be presented by a finite sequence of {\em detour moves}, described in Fig.~\ref{fig:detour}. Each detour move can be presented by a finite sequence of moves (1'), (2'), (3'), (3''), (4') in Fig. \ref{fig:moves2}.
    \begin{figure}[h]
\begin{center}
 \includegraphics[width = 6cm]{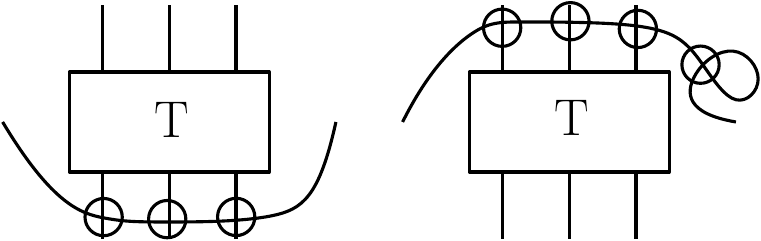}

\end{center}

\caption{A detour move: an arc only with virtual crossings can pass through classical crossings
and double lines.}\label{fig:detour}
\end{figure}
\end{proof}

\subsection{Braiding process}
Let $O$ be the origin of $\mathbb{R}^{2}$. Then $\mathbb{R}^{2}-\{O\}$ has the open book decomposition $\mathbb{R}_{+} \times S^{1}$ by polar coordinate, where $\mathbb{R}_{+}$ is the set of positive real numbers. Let $\pi: \mathbb{R}_{+} \times S^{1} \rightarrow S^{1}$ denote the second coordinate projection.

For a diagram $D$ with double lines we consider $\dL(D)$ as the set of points where double lines are placed. 

\begin{dfn}
    A diagram $D$ with double lines is called {\em a closed dl-braid diagram} if it satisfies the following conditions:
    \begin{itemize}
        \item[(1)] It is contained in $\mathbb{R}^{2}-\{O\} \cong \mathbb{R}_{+} \times S^{1}$.
        \item[(2)] Let $K : \sqcup S^{1} \rightarrow \mathbb{R}^{2}-\{O\}$ be the underlying immersion of $D$. Then $\pi\circ K$ is a covering map of $S^{1}$ of degree $n$ with respect to the orientations of $\sqcup S^{1}$ and $S^{1}$.
    \end{itemize}  
    A closed dl-braid diagram is {\em good} if it satisfies the following condition:
    \begin{itemize}
        \item[(3)]  Let $N_{1},\dots, N_{m}$ be regular neighborhoods of the classical crossings and double lines of $D$. Then $\pi(N_{i}) \cap \pi(N_{j}) = \emptyset$ for $i\neq j$.
    \end{itemize}
\end{dfn}

\begin{prop}\label{prop:braiding-process}
    Every diagram $D$ with double lines is equivalent to a good closed dl-braid diagram $D'$ such that $D$ and $D'$ have the same Gauss data.
\end{prop}

\begin{proof}
    The proof is analogous to the proof of Proposition 3.5 in \cite{NegiPrabhakarKamada}. The idea is, first we move $D$ so that $D$ is placed in $\mathbb{R}^{2}-\{O\}$ and $N_{1},\dots,N_{m}$ satisfy that $\pi(N_{i}) \cap \pi(N_{j}) = \emptyset$ for $i\neq j$. Second, we move $D$ so that arcs in $N_{i}$ agrees with the orientation of $S^{1}$ of $\mathbb{R}^{2}-\{O\} \cong \mathbb{R}_{+} \times S^{1}$. Then we replace arcs of $D$ in $W(D)$ so that the result is a good closed dl-braid diagram. It is possible since the detour move can be applied.
\end{proof}

We call the procedure in the proof of Proposition~\ref{prop:braiding-process} {\em the braiding process.}
\begin{thm}
    Any link is represented as the closure of a virtual dl-braid diagram.
\end{thm}

\subsection{Markov theorem for virtual dl-braids}
{\em A dl-Markov move of type $0$} or {\em a dl-M0-move} is a move for virtual dl-braids described in Fig.~\ref{fig:dl-braid-moves}.

{\em A dl-Markov move of type $1$} or {\em a dl-M1-move} is a conjugation of a virtual dl-braid by another virtual dl-braid, that is, $\beta \sim \beta_{1}^{-1}\beta \beta_{1}$.  

Let us define a map $l_{s}^{t}(\beta): VB_{n}^{dl} \rightarrow VB_{n+s+t}^{dl}$ by adding $s$ strands on the left of $\beta$ and $t$ strands on the right of $\beta$.
{\em A dl-Markov move of type $2$} or {\em a dl-M2-move} is stabilizations of positive type $\beta \sim l_{0}^{1}(\beta)\sigma_{n}$, negative type $\beta \sim l_{0}^{1}(\beta)\sigma_{n}^{1}$, or virtual type $\beta \sim l_{0}^{1}(\beta)\rho_{n}$ or their inverse operations.

{\em A dl-Markov move of type $3$} or {\em a dl-M3-move} is the left and right virtual exchange moves, that is, $$l_{1}^{0}(\beta_{1})\sigma_{1}^{-1} l_{1}^{0}(\beta_{2})\sigma_{1} \sim l_{1}^{0}(\beta_{1})\rho_{1}^{-1} l_{1}^{0}(\beta_{2})\rho_{1},$$

and $$l_{0}^{1}(\beta_{1})\sigma_{n}^{-1} l_{0}^{1}(\beta_{2})\sigma_{n} \sim l_{0}^{1}(\beta_{1})\rho_{n}^{-1} l_{0}^{1}(\beta_{2})\rho_{n}.$$

\begin{dfn}
    Two virtual dl-braid diagrams $\beta_{1}$ and $\beta_{2}$ are {\em dl-Markov equivalent} if one can be obtained from another by a finite sequence of dl-M0-, dl-M1-, dl-M2- and dl-M3-moves.
\end{dfn}

\begin{thm}\label{thm:dl-Markov}
    Two virtual dl-braids have equivalent closures as knots in $S_{g}\times S^{1}$ if and only if they are dl-Markov equivalent.
\end{thm}

The above notions can be rewritten in the point of view of closed virtual dl-braids. Let $D$ and $D'$ be closed virtual dl-braid diagrams and $\beta$ and $\beta'$ be virtual dl-braid diagrams obtained by cutting $D$ and $D'$ along $\pi^{-1}(\theta)$ and $\pi^{-1}(\theta')$ for some $\theta,\theta' \in [0,2\pi)$. Let us say that $D'$ is obtained from $D$ by {\em a dl-Markov move of type $0$} or {\em a dl-M0-move} if $\beta'$ is obtained from $\beta$ by a finite sequence of dl-M0-move and dl-M1-move. Let us say that $D'$ is obtained from $D$ by {\em a dl-Markov move of type $2$} or {\em a dl-M2-move} if $\beta'$ is obtained from $\beta$ by a finite sequence of dl-M2-move and dl-M1-move. And we say that $D'$ is obtained from $D$ by {\em a dl-Markov move of type $3$} or {\em a dl-M3-move} if $\beta'$ is obtained from $\beta$ by a finite sequence of dl-M3-move and dl-M1-move.

\begin{dfn}
    Two closed virtual dl-braid diagrams $D$ and $D'$ are {\em dl-Markov equivalent} if one can be obtained from another by a finite sequence of dl-M0-, dl-M2-, dl-M3-moves.
\end{dfn}

\begin{lem}\label{lem:cl-dl-Markov}
    Two closed virtual dl-braid diagrams $D$ and $D'$ are dl-Markov equivalent if and only if $\beta$ and $\beta'$ are dl-Markov equivalent, where $\beta$ and $\beta'$ are virtual dl-braid diagrams obtained by cutting $D$ and $D'$ along $\pi^{-1}(\theta)$ and $\pi^{-1}(\theta')$.
\end{lem}

\begin{proof}
    It is clear by definitions.
\end{proof}

From Lemma~\ref{lem:cl-dl-Markov} it follows that Theorem~\ref{thm:dl-Markov} is equivalent to the following theorem.
\begin{thm}\label{thm:cl-dl-Markov}
    Two closed virtual dl-braid diagrams are equivalent as knots in $S_{g}\times S^{1}$ if and only if they are dl-Markov equivalent.
\end{thm}

To prove Theorem~\ref{thm:cl-dl-Markov}, the following Lemma is required.

\begin{lem}\label{lem:dlbraid-GD-dlMarkov}
    Two closed virtual dl-braid diagrams with the same Gauss data are dl-Markov equivalent.
\end{lem}
\begin{proof}
    The proof is exactly same with the proof of Lemma 4.7 in \cite{NegiPrabhakarKamada}. The idea is as follows: Assume that $D$ and $D'$ have the same Gauss data. Then there exists a one-to-one correspondence between regular neighborhoods $\{N(v_{1}),\dots, N(v_{m})\}$ of crossings and double lines of $D$ and regular neighborhoods $\{N(v'_{1}),\dots, N(v'_{m})\}$ of crossings and double lines of $D'$. By dl-M0- and dl-M2-moves one can deform them so that $\{N(v_{1}),\dots, N(v_{m})\}$ and $\{N(v'_{1}),\dots, N(v'_{m})\}$ are placed in the same cyclic order. Notice that in this step we mainly use stabilization of virtual type and detour moves. So, our special moves (4), (4') and (5) related with double lines do not affect on this process. And then we deform the arcs on $W(D)$ and $W(D')$ so that they have the same form by using dl-M0- and dl-M2-moves. Again, we mainly use detour moves.
\end{proof}

\begin{proof}[Proof of Theorem~\ref{thm:cl-dl-Markov}]
    It is clear that if two closed virtual dl-braid diagrams are dl-Markov equivalent, then they are equivalent as knots in $S_{g}\times S^{1}$. Conversely, suppose that $D$ and $D'$ are closed virtual dl-braid diagrams which are equivalent as knots in $S_{g}\times S^{1}$. There exists a finite sequence of diagrams with double lines 
    $$D=D_{0} \rightarrow D_{1} \rightarrow \cdots \rightarrow D_{m} = D',$$
    such that $D_{i+1}$ is obtained from $D_{i}$ by one of moves in Fig.~\ref{fig:moves2}. Let $\tilde{D}_{i}$ be a closed virtual dl-braid diagram obtained from $D_{i}$ by the braiding process. Since $D_{i}$ and $\tilde{D}_{i}$ have the same Gauss data, it is sufficient to show that $\tilde{D}_{i}$ and  $\tilde{D}_{i+1}$ are dl-Markov equivalent. It is shown in \cite{Kamada-braid} that if $D_{i+1}$ is obtained from $D_{i}$ by moves (1),(2),(3),(4),(1'),(2'),(3') and (3'') in Fig.~\ref{fig:moves2}, that is, generalized Reidemeister moves, then $\tilde{D}_{i}$ and $\tilde{D}_{i+1}$ are dl-Markov equivalent. For the moves (5) and (4') in Fig.~\ref{fig:moves2} the situations are same with the moves TM1 and TM2 for twisted virtual links. Let us assume that $D_{i+1}$ is obtained from $D_{i}$ by the move (4). Let $\{N(v_{1}),\dots, N(v_{l})\}$ be regular neighborhoods of crossings and double lines of $\tilde{D}_{i}$ and $\{N(v'_{1}),\dots, N(v'_{l})\}$ of crossings and double lines of $\tilde{D}_{i+1}$. Let us say $v_{1}$ and $v_{2}$ be the crossing and the double lines of $\tilde{D}_{i}$ contained in move (4) and let $v_{1}'$ and $v_{2}'$ be the crossing and the double lines of $\tilde{D}_{i+1}$ contained in move (4). Since $N(v_{1})$ and $N(v_{2})$ are placed closely, let us take a neighborhood $\hat{N}(v_{1,2}) = N(v_{1}) \cup N(v_{2}) \cup N(\alpha)$, where $N(\alpha)$ is a regular neighborhood of the arc $\alpha$ connecting $v_{1}$ and $v_{2}$. We can apply the braiding process to $D_{i}$ so that $\{\hat{N}(v_{1,2}), N(v_{3}),\dots, N(v_{l})\}$ are placed in cyclic order and obtain new closed virtual dl-braid diagram $\tilde{\tilde{D}}_{i}$. Since $\tilde{D}_{i}$ and $\tilde{\tilde{D}}_{i}$ have the same Gauss data, they are Markov equivalent by Lemma~\ref{lem:dlbraid-GD-dlMarkov}. Similarly, one can obtain $\tilde{\tilde{D}}_{i+1}$ so that $\{\hat{N}(v_{1,2}'), N(v_{3}'),\dots, N(v_{l}')\}$ are placed in this cyclic order. Then $\tilde{D}_{i+1}$ and $\tilde{\tilde{D}}_{i+1}$ have the same Gauss data and hence they are dl-Markov equivalent. Since $D_{i}$ and $D_{i+1}$ have the same diagrams except $\hat{N}(v_{1,2})$ and $\hat{N}(v_{1,2}')$, we can assume that $\tilde{\tilde{D}}_{i}$ and $\tilde{\tilde{D}}_{i+1}$ have the same diagram except $\hat{N}(v_{1,2})$ and $\hat{N}(v_{1,2}')$. Since $\hat{N}(v_{1,2}') \cap \tilde{\tilde{D}}_{i+1}$ can be obtained by applying the move (4) to $\hat{N}(v_{1,2}) \cap \tilde{\tilde{D}}_{i}$, $\tilde{\tilde{D}}_{i}$ and $\tilde{\tilde{D}}_{i+1}$ are dl-Markov equivalent and hence $\tilde{D}_{i}$ and $\tilde{D}_{i+1}$ are dl-Markov equivalent.
\end{proof}

\end{document}